\documentclass[12pt,twoside]{amsart}
\usepackage{mathrsfs}
\usepackage{stmaryrd}
\usepackage{amsmath}
\usepackage{amsthm}
\usepackage{amsfonts}
\usepackage{amssymb}
\usepackage{latexsym}
\usepackage[all]{xy}
\usepackage{epsfig}
\usepackage{amscd,mathrsfs,graphicx,mathrsfs}

\date{}
\allowdisplaybreaks[4] \footskip=15pt

\renewcommand{\uppercasenonmath}[1]{}

\numberwithin{equation}{section} \theoremstyle{plain}
\newtheorem*{thm*}{Main Theorem}
\newtheorem{thm}{Theorem}[section]
\newtheorem{cor}[thm]{Corollary}
\newtheorem*{cor*}{Corollary}
\newtheorem{lem}[thm]{Lemma}
\newtheorem*{lem*}{Lemma}
\newtheorem{prop}[thm]{Proposition}
\newtheorem*{prop*}{Proposition}
\newtheorem{fact}[thm]{Fact}
\newtheorem*{fact*}{Fact}
\newtheorem{rem}[thm]{Remark}
\newtheorem*{rem*}{Remark}
\newtheorem{exa}[thm]{Example}
\newtheorem*{exa*}{Example}
\newtheorem{df}[thm]{Definition}
\newtheorem*{df*}{Definition}
\newtheorem*{ack*}{ACKNOWLEDGEMENTS}

\newcommand{\Ext}{\mbox{\rm Ext}}
\newcommand{\Hom}{\mbox{\rm Hom}}

\newcommand{\pd}{\mbox{\rm pd}}

\newcommand{\Ker}{\mbox{\rm Ker}}
\newcommand{\res}{\mbox{\rm res}}

\begin{document}
\begin{center}
{\Large \bf RELATIVE DERIVED CATEGORY WITH RESPECT TO A SUBCATEGORY}
 \footnotetext{* Corresponding author.

E-mail:~dizhenxing19841111@126.com,~z990303@seu.edu.cn,~renwei@nwnu.enu.cn,~and\\jlchen@seu.edu.cn}

\vspace{0.5cm}  ZHENXING DI$^{1}$, XIAOXIANG ZHANG$^{1}$, WEI REN$^{2}$ and JIANLONG CHEN$^{1,\,*}$ \\\

1.  \emph{Department of Mathematics, Southeast University}, \\\emph{Nanjing} \rm210096, \emph{P. R. China} \\

2.  \emph{Department of Mathematics}, \emph{Northwest Normal University}, \\\emph{Lanzhou} \rm730070, \emph{P. R. China} \\
\end{center}


\bigskip
\centerline { \bf  Abstract}
The notion of relative derived category with respect to a subcategory is introduced.
A triangle-equivalence, which extends a theorem of Gao and Zhang
[Gorenstein derived categories, \emph{J. Algebra} \textbf{323} (2010) 2041-2057] to the bounded below case, is obtained.
Moreover, we interpret the relative derived functor $\mathrm{Ext}_{\mathcal{X}\mathcal {A}}(-,-)$ as the morphisms in such derived category
and give two applications.
\bigskip
\leftskip10truemm \rightskip10truemm \noindent   \\
\vbox to 0.3cm{}

\leftskip0truemm \rightskip0truemm \vbox to 0.2cm{}

\section { \bf Introduction}

For an abelian category $\mathcal {A}$ with enough projective objects,
the Gorenstein (projective) derived category,
which makes Gorenstein projective quasi-isomorphisms become isomorphisms,
was introduced by Gao and Zhang \cite{7}
to close a gap of the corresponding version of derived category in so-called Gorenstein homological algebra.
In particular, they established the triangle-equivalence
$$\indent\indent\indent\mathrm{\mathbf{K}}^{b}(\mathcal {G}(\mathcal {P}))\cong
\mathrm{\mathbf{D}}^{b}_{\mathcal {G}(\mathcal {P})}(\res\,\widehat{\mathcal {G}(\mathcal {P})})\indent\indent\indent(\dag)$$
where $\res\,\widehat{\mathcal {G}(\mathcal {P})}$ denotes the full subcategory of objects in $\mathcal {A}$
with finite Gorenstein projective dimension (see \cite[Theorem 3.6]{7}).
Moreover, they showed, as one of the advantages of such relative derived category,
the relative derived functor Ext$_{\mathcal {G}(\mathcal {P})}(-,-)$,
which is derived from Hom$_\mathcal {A}(-,-)$ using proper $\mathcal {G}(\mathcal {P})$-resolutions of the first variable,
can be interpreted as the morphisms in $\mathrm{\mathbf{D}}^{b}_{\mathcal {G}(\mathcal {P})}(\mathcal {A})$
(see \cite[Theorem 3.12]{7}).

Let $\mathcal {X}$ and $\mathcal {S}$ be subcategories of $\mathcal {A}$
with $\mathcal {S}$ closed under direct summands.
We introduce in the paper the notion of relative derived category with respect to $\mathcal {X}$,
$\mathrm{\mathbf{D}}^{\ast}_{\mathcal{X}}(\mathcal {S})$ for $\ast \in \{{blank, -, b}\}$,
which is defined to be the Verdier quotient of the homotopy category $\mathrm{\mathbf{K}}^{\ast}(\mathcal {S})$ with respect to
the thick triangulated subcategory of all $\mathcal {X}$-acyclic complexes in $\mathrm{\mathbf{K}}^{\ast}(\mathcal {S})$
(see Definition \ref{def 3.4}).
Then we obtain in section 3 the following triangle-equivalence,
which extends $(\dag)$ to the bounded below case
when we specially take $\mathcal {X}=\mathcal {G}(\mathcal {P})$ (see Theorem \ref{thm 3.10}).
Also, it can be compared to \cite[Proposition 3.5]{4} (see Example \ref{exa 3.10}).\vspace{2mm}\\
$\mathbf{Theorem}\,\mathbb{I}\,$ \emph{Let $\mathcal {X}$ be a subcategory of $\mathcal {A}$.
Assume that $\mathcal {X}$ is exact and has an injective cogenerator.
Then there exists a triangle-equivalence
$\mathrm{\mathbf{K}}^{-}(\mathcal {X})\cong \mathrm{\mathbf{D}}^{-}_{\mathcal{X}}(\res\,\widehat{\mathcal {X}})$}.\vspace{2mm}

The study of relative homological algebra goes back to Butler and Horrocks \cite{B}, and Eilenberg and Moore \cite{EM}.
It was reinvigorated recently by a number of authors (see, for example, Enochs and Jenda \cite{EJ2000}, and Avramov and Martsinkovsky \cite{1}).
Assume that $M$ and $N$ are objects in $\mathcal {A}$ such that $M$ admits a proper $\mathcal {X}$-resolution $\mathbf{X}$.
According to \cite{16}, the relative cohomology groups $\mathrm{Ext}^{n}_{\mathcal{X}\mathcal {A}}(M,N)$
is defined as $\mathrm{H}^{n}\,(\Hom_{\mathcal {A}}(\mathbf{X},N))$.
In section 4, inspired by \cite[Theorem 3.12]{7},
we give the following new characterization of $\mathrm{Ext}^{n}_{\mathcal{X}\mathcal {A}}(M, N)$ (see Theorem \ref{thm 4.2}).\vspace{2mm}\\
$\mathbf{Theorem}\,\mathbb{II} \,$ \emph{Let $\mathcal {X}$ be a subcategory of $\mathcal {A}$ and $M$, $N$ objects in $\mathcal {A}$.
Assume that $M$ admits a proper $\mathcal {X}$-resolution.
Then $\mathrm{Ext}^{n}_{\mathcal{X}\mathcal {A}}(M, N) = \Hom_{\mathrm{\mathbf{D}}^{b}_{\mathcal{X}}(\mathcal {A})}(M, N[n])$.}\vspace{2mm}

Moreover, as an application of Theorem $\mathbb{II}$,
we give in derived category a new proof for the existence of the Avramov-Martsinkovsky type exact sequence appeared in \cite[Theorem B]{18}
(see Corollary \ref{cor 4.7}).

In the next section, we mainly review some definitions and notation.
It should be pointed out that many of them were used by Sather-Wagstaff et al. [19-22].

\section { \bf Preliminaries}

Throughout this paper, $\mathcal {A}$ denotes an abelian category with enough projective objects.
The term \emph{subcategory} stands for a full additive subcategory of $\mathcal {A}$ closed under isomorphisms.
Following \cite{17}, a subcategory of $\mathcal {A}$ is called \emph{exact} if it is closed under extensions and direct summands.
We write $\mathcal {P}=\mathcal {P}(\mathcal {A})$ for the subcategory of projective objects in $\mathcal {A}$.

Let $\mathcal {X},\,\mathcal {Y}$ and $\mathcal {W}$ be subcategories of $\mathcal {A}$
such that $\mathcal {W}\subseteq\mathcal {X}$.
We write $\mathcal {X}~\bot~\mathcal {Y}$ in case Ext$^{\geqslant1}_\mathcal {A}(X,Y)=0$ for each object $X\in\mathcal {X} $ and
each object $Y\in\mathcal {Y}$.
When $\mathcal{X} = \{M\}$,
we use the notation $M \,\bot\, \mathcal{Y}$ instead of $\{M\} \,\bot\, \mathcal{Y}$.
There is an analogue $\mathcal{X} \,\bot\, M$.
According to \cite{17}, $\mathcal {W}$ is said to be a \emph{cogenerator} for $\mathcal {X}$
if for each object $X\in\mathcal {X}$,
there exists a short exact sequence
$0\to X\to W\to X'\to0$
such that $W\in\mathcal {W}$ and $X'\in\mathcal {X}$.
$\mathcal{W}$ is called an \emph{injective cogenerator} for $\mathcal {X}$
if $\mathcal {W}$ is a cogenerator for $\mathcal {X}$ and $\mathcal{X}~\bot~\mathcal {W}$.
\emph{Generator} and \emph{projective generator} can be defined dually.

A \emph{complex} $\mathbf{X}$ is often displayed as a sequence of objects in $\mathcal {A}$
$$\xymatrix{
  \cdots\ar[r] & X^{-1} \ar[r]^{\delta_\mathbf{X}^{-1}} & X^{0} \ar[r]^{\delta_\mathbf{X}^{0}}
  & X^{1} \ar[r]& \cdots} $$
with $\delta_\mathbf{X}^{n+1}\delta_\mathbf{X}^{n} = 0$ for all $n\in \mathbb{Z}$.
The $n$th \emph{homology} of the complex $\mathbf{X}$ is defined as
$\mathrm{H}^{n}(\mathbf{X}) = \mathrm{Ker}(\delta_\mathbf{X}^{n}) / \mathrm{Im}(\delta_\mathbf{X}^{n-1})$.
We identify $M$ with the complex
$\cdots \to 0\to M\to 0\to \cdots$,
where $M$ is in degree zero and 0 elsewhere.
For an integer $n$, $\mathbf{X}[n]$ denotes the complex $\mathbf{X}$ shifting $n$ degree, that is,
$X[n]^{m} = X^{n+m}$ and $\delta_{\mathbf{X}[n]}^{m} = (-1)^{n}\delta^{n+m}_{\mathbf{X}}$.
Given two complexes $\mathbf{X}$ and $\mathbf{Y}$,
the complex $\Hom_\mathcal {A}(\mathbf{X}, \mathbf{Y})$ is defined with
$\Hom_\mathcal {A}(\mathbf{X}, \mathbf{Y})^{n} = \prod_{k\in \mathbb{Z}}\Hom_\mathcal {A}(X^{k}, Y^{k+n})$,
and with differential $\delta^{n}((f^{k})_{k\in \mathbb{Z}}) = (\delta^{k+n}_{ \mathbf{Y}}f^{k} - (-1)^{n}f^{k+1}\delta^{k}_{\mathbf{X}})_{k\in\mathbb{Z}}$
for $(f^{k})_{k\in \mathbb{Z}}\in \Hom_\mathcal {A}(\mathbf{X}, \mathbf{Y})^{n}$.

A \emph{morphism} $f: \mathbf{X}\to \mathbf{Y}$ of complexes is a family of morphisms
$f = (f^{n}: X^{n}\to Y^{n})_{n\in \mathbb{Z}}$ of objects in $\mathcal {A}$
satisfying $\delta_{\mathbf{Y}}^{n}f^{n} = f^{n+1}\delta^{n}_{\mathbf{X}}$
for all $n\in \mathbb{Z}$.
Morphisms $f,\, g: \mathbf{X}\to \mathbf{Y}$ are called \emph{homotopic}, denoted by $f\sim g$,
if there exists a family of morphisms $(s^{n}: X^{n}\to Y^{n-1})_{n\in\mathbb{Z}}$ of objects in $\mathcal {A}$
satisfying $f^{n}-g^{n} = \delta^{n-1}_{\mathbf{Y}}s^{n} + s^{n+1}\delta^{n}_{\mathbf{X}}$
for all $n\in \mathbb{Z}$.
A \emph{quasi-isomorphism} is a morphism $f: \mathbf{X}\to \mathbf{Y}$ with
$\textrm{H}^{n}(f): \textrm{H}^{n}(\mathbf{X})\to \textrm{H}^{n}(\mathbf{Y})$
bijective for all $n\in \mathbb{Z}$.

A complex $\mathbf{X}$ is called an $\mathcal{X}$-\emph{resolution} of $M$
if $X^i \in \mathcal{X}$ for all $i \leqslant 0$, $X^i = 0$ for all $i > 0$,
H$^i\,(\mathbf{X})$ = 0 for all $i < 0$ and H$^0\,(\mathbf{X}) \cong M$.
In this case, the associated exact sequence
$$\cdots \to X^{-1} \to X^0 \to M \to 0$$
is denoted by \textbf{X}$^+$.
Sometimes we call the quasi-isomorphism
$\xymatrix@C=0.5cm{
    \mathbf{X} \ar[r]^{\simeq} & M}$
an $\mathcal {X}$-resolution of $M$.
If \textbf{X}$^+$ is Hom$_\mathcal {A}(\mathcal{X}, -)$-exact,
we say that \textbf{X} is a \emph{proper} $\mathcal{X}$-resolution of $M$,
and let $\res\,\widetilde{\mathcal {X}}$ denote the subcategory of objects in $\mathcal {A}$ admitting a proper $\mathcal {X}$-resolution.

The $\mathcal{X}$-\emph{projective dimension} of $M$ is the least non-negative $n$ such that there exists an exact sequence
$0 \to X^{-n} \to \cdots \to X^{-1} \to X^0 \to M \to 0$,
where each $X^{-i} \in \mathcal{X}$.
In this case, we write $\mathcal{X}$-pd$(M) = n$.
If no such $n$ exists, we write $\mathcal{X}$-pd$(M) = \infty$.
The subcategory of objects in $\mathcal {A}$ with finite $\mathcal{X}$-projective dimension is denoted by $\res\,\widehat{\mathcal{X}}$.

For $\ast \in \{{blank, -, b}\}$,
$\mathrm{\mathbf{K}}^{\ast}(\mathcal {A})$ and $\mathrm{\mathbf{D}}^{\ast}(\mathcal {A})$
stand for the corresponding homotopy category and derived category of $\mathcal {A}$, respectively.
We will often use the standard formula
$\Hom_{\mathrm{\mathbf{K}}^{\ast}(\mathcal {A})}(\mathbf{X}, \mathbf{Y}[n]) = \textrm{H}^{n}\, (\Hom_{\mathcal {A}}(\mathbf{X}, \mathbf{Y}))
= \Hom_{\mathrm{\mathbf{K}}^{\ast}(\mathcal {A})}(\mathbf{X}[-n], \mathbf{Y})$.

Let $\mathcal{B}$ be a triangulated subcategory of a triangulated category $\mathcal{K}$,
and let $\Sigma$ be the compatible multiplicative system determined by $\mathcal{B}$.
In the Verdier quotient  \cite[Chapter 2]{13} $\mathcal{K} / \mathcal{B} = \Sigma^{-1}\mathcal{K}$,
each morphism $f: X\to Y$ is given by an equivalence class of right fractions $a/s$
presented by $X\stackrel{s}\Longleftarrow Z \stackrel{a}\longrightarrow Y$,
where the double arrow means $s\in \Sigma$.

\section{\bf Relative Derived Categories}

In what follows, we always assume that $\mathcal {X}$ and $\mathcal {S}$ are subcategories of $\mathcal {A}$
with $\mathcal {S}$ closed under direct summands.
We begin with the following definition.

\begin{df}\label{def 3.1}
\rm{A complex $\mathbf{S}$ is called $\mathcal{X}$-\emph{acyclic}
 if $\Hom_{\mathcal {A}}( X,\mathbf{S})$ is acyclic for each object $X\in\mathcal {X}$.
 A morphism $f: \mathbf{S}\to \mathbf{T}$ of complexes is called an $\mathcal{X}$\emph{-quasi-isomorphism}
 if $\Hom_{\mathcal {A}}(X, f)$ is a quasi-isomorphism for each object $X\in\mathcal {X}$.}
 \end{df}

Let $\mathcal {T}$ be a subcategory of $\mathcal {A}$.
For $\ast \in \{{blank, -, b}\}$, we denote by $\mathrm{\mathbf{K}}^{\ast}(\mathcal{T})$
the homotopy category with each complex constructed by objects in $\mathcal{T}$
and by $\mathrm{\mathbf{K}}^{\ast}_{\mathcal{X}}(\mathcal{T})$
the subcategory of $\mathrm{\mathbf{K}}^{\ast}(\mathcal{T})$
consisting of all $\mathcal{X}$-acyclic complexes.

\begin{fact}\label{rem 3.2}
         By virtue of \cite[Lemma 2.4]{3}, a complex $\mathbf{S}$ is $\mathcal{X}$-acyclic
         if and only if $\Hom_{\mathcal {A}}(\mathbf{D}, \mathbf{S})$ is acyclic for each complex $\mathbf{D}\in \mathrm{\mathbf{K}}^{-}(\mathcal {X})$.
         Moreover, it follows from \cite[Proposition 2.6]{3} that
         a morphism $f: \mathbf{S}\to \mathbf{T}$ of complexes is an $\mathcal{X}$-quasi-isomorphism
         if and only if $\Hom_{\mathcal {A}}(\mathbf{D}, f)$ is a quasi-isomorphism
         for each complex $\mathbf{D}\in \mathrm{\mathbf{K}}^{-}(\mathcal{X})$.
\end{fact}

Recall from \cite{8} that a full triangulated subcategory $\mathcal{C}$ of a triangulated  category $\mathcal{D}$ is said to be \emph{thick}
if it satisfies the following condition:
assume that a morphism $f: X\to Y$ in $\mathcal{D}$ can be factored through an object from $\mathcal{C}$,
and enters a distinguished triangle $X\stackrel{f}\to Y\to Z\to X[1]$ with $Z$ in $\mathcal{C}$,
then $X$ and $Y$ are objects in $\mathcal{C}$.
A standard example of thick subcategory is the subcategory of all acyclic complexes in $\mathrm{\mathbf{K}}(\mathcal {A})$,
that is, $\mathrm{\mathbf{K}}^{\ast}_{\mathcal{P}}(\mathcal {A})$.

The following result depends on an important characterization of thick subcategories due to Rickard, called \emph{Rickard's criterion}:
a full triangulated subcategory $\mathcal{C}$ of a triangulated category $\mathcal{D}$ is thick
if and only if every direct summand of an object of $\mathcal{C}$ is in $\mathcal{C}$
(see \cite[Proposition 1.3]{14} or \cite[Criterion 1.3]{12}).

\begin{lem}\label{lem 3.3}
For $\ast \in \{{blank, -, b}\}$, $\mathrm{\mathbf{K}}^{\ast}_{\mathcal{X}}(\mathcal {S})$
is a thick subcategory of $\mathrm{\mathbf{K}}^{\ast}(\mathcal {S})$.
\end{lem}

\begin{proof}
By Rickard's criterion, it suffices to show that
$\mathrm{\mathbf{K}}^{\ast}_{\mathcal{X}}(\mathcal {S})$
is a full triangulated subcategory of $\mathrm{\mathbf{K}}^{\ast}(\mathcal {S})$ and it is closed under direct summands.
This is clear.
\end{proof}

Note that a morphism $f: \mathbf{S}\to \mathbf{T}$ of complexes is an $\mathcal{X}$-quasi-isomorphism
if and only if its mapping cone cone$(f)$ is $\mathcal{X}$-acyclic.
Then, the collection of all $\mathcal{X}$-quasi-isomorphisms in $\mathrm{\mathbf{K}}^{\ast}(\mathcal {S})$,
denoted by $\Sigma^{\mathcal {S}}_{\mathcal{X}}$,
is a saturated multiplicative system corresponding to the subcategory
$\mathrm{\mathbf{K}}^{\ast}_{\mathcal{X}}(\mathcal {S})$.

\begin{df}\label{def 3.4}
{\rm
For $\ast \in \{{blank, -, b}\}$, the \emph{relative derived category} $\mathrm{\mathbf{D}}^{\ast}_{\mathcal{X}}(\mathcal {S})$
is defined to be the Verdier quotient of $\mathrm{\mathbf{K}}^{\ast}(\mathcal {S})$,
that is,
$$\mathrm{\mathbf{D}}^{\ast}_{\mathcal{X}}(\mathcal {S}):
= \mathrm{\mathbf{K}}^{\ast}(\mathcal {S}) / \mathrm{\mathbf{K}}^{\ast}_{\mathcal{X}}(\mathcal {S})
= \Sigma^{-\mathcal {S}}_{\mathcal{X}}\, \mathrm{\mathbf{K}}^{\ast}(\mathcal {S}).$$

Note that $\mathrm{\mathbf{D}}^{\ast}_{\mathcal{X}}(\mathcal {S})$ is the derived category in sense of Neeman \cite{12},
of the exact category $(\mathcal {S}$, $\mathcal{E}^{\mathcal {S}}_{\mathcal{X}})$,
where $\mathcal{E}^{\mathcal {S}}_{\mathcal{X}}$ consists of all short $\mathcal{X}$-acyclic sequences in $\mathcal {S}$.}
\end{df}

The next result will be used in the proof of Lemma \ref{lem 3.7}.

\begin{lem}\label{lem 3.6}
Let $\mathbf{S}$ be a complex.
Assume that there is a complex $\mathbf{D}\in \mathrm{\mathbf{K}}^{-}(\mathcal{X})$
and an $\mathcal{X}$-quasi-isomorphism $f: \mathbf{S}\to \mathbf{D}$.
Then there exists a morphism $g: \mathbf{D}\to \mathbf{S}$ such that $fg$ is homotopic to ${\rm Id}_{\mathbf{D}}$.
\end{lem}

\begin{proof}
In view of Fact \ref {rem 3.2},
$\Hom_{\mathcal {A}}(\mathbf{D}, f): \Hom_{\mathcal {A}}(\mathbf{D}, \mathbf{S}) \to \Hom_{\mathcal {A}}(\mathbf{D},\mathbf{D})$
is a quasi-isomorphism.
Then by \cite [(1.1)]{1}, for the morphism ${\rm Id}_{\mathbf{D}}$,
there is a morphism $g: \mathbf{D}\to \mathbf{S}$ such that $fg \sim {\rm Id}_{\mathbf{D}}$ .
\end{proof}

By a slight modification of the proof of \cite[Proposition 2.8]{7},
we can get the following result, which makes the morphisms in $\mathrm{\mathbf{D}}^-_{\mathcal{X}}(\mathcal {S})$ easier to understand.

\begin{lem}\label{lem 3.7}
Let $\mathbf{D}$ be a complex in $\mathrm{\mathbf{K}}^{-}(\mathcal{X})$ and $\mathbf{S}$ a complex in $\mathrm{\mathbf{K}}^-(\mathcal{S})$.
Assume $\mathcal {X}\subseteq\mathcal {S}$.
Then $\varphi: f\to f/{\rm Id}_{\mathbf{D}}$ gives an isomorphism of abelian groups
$\Hom_{\mathrm{\mathbf{K}}^-(S)}(\mathbf{D}, \mathbf{S})\cong \Hom_{\mathrm{\mathbf{D}}^-_{\mathcal{X}}(\mathcal {S})}(\mathbf{D}, \mathbf{S})$.
\end{lem}
\begin{proof}
If $f/{\rm Id}_{\mathbf{D}} =0$, then by the calculus of right fractions there is an $\mathcal{X}$-quasi-isomorphism
$s: \mathbf{Y}\to \mathbf{D}$
for some complex $\mathbf{Y}$ such that $fs \sim 0$.
It follows from Lemma \ref{lem 3.6} that there is a morphism $g: \mathbf{D}\to \mathbf{Y}$ such that $sg \sim {\rm Id}_{\mathbf{D}}$.
Thus, $f\sim fsg \sim 0$.
Moreover, for each $f/s \in \Hom_{\mathrm{\mathbf{D}}^-_{\mathcal{X}}(\mathcal {S})}(\mathbf{D}, \mathbf{S})$
presented by $\mathbf{D}\stackrel{s}\Longleftarrow \mathbf{Y} \stackrel{f}\longrightarrow \mathbf{S}$,
by Lemma \ref{lem 3.6} there is a morphism $g: \mathbf{D}\to \mathbf{Y}$ such that $sg \sim {\rm Id}_{\mathbf{D}}$.
This implies that $f/s = fg/{\rm Id}_{\mathbf{D}} = \varphi (fg)$.
Thus, $\varphi$ is an isomorphism, as desired.
\end{proof}

The next result is given in service of Lemma \ref{lem 3.9}.

\begin{lem}\label{lem 3.8}
Let $\mathbb{L}=\,\, 0\to
\Ker(g)\to\mathbf{M}\overset{g
}{\longrightarrow}\mathbf{N}\to 0$ be a short exact
sequence of complexes such that $\Ker(g)$ is $\mathcal
{X}$-acyclic and $\mathbb{L}$ remains exact after applying the
functor $\Hom_\mathcal {A}(\mathcal {X},-)$, i.e., the sequence
$$\indent\indent 0\to\Hom_\mathcal {A}(X,\Ker(g))\to \Hom_\mathcal {A}(X,\mathbf{M})\to
\Hom_\mathcal {A}(X,\mathbf{N})\to 0
\indent\indent(\ddag)$$ of complexes is exact for each object
$X\in\mathcal {X}$. Then for each complex $\mathbf{D}\in
\mathrm{\mathbf{K}}^{-}(\mathcal {X})$ and each morphism
$\alpha:\mathbf{D} \to \mathbf{N}$, there exists a
morphism $\beta:\mathbf{D} \to \mathbf{M}$ such that the
following diagram
$$
\xymatrix{ &&& \mathbf{D}\ar[d]_{\alpha}\ar@.[dl]_{\beta}\\
          0\ar[r]& \Ker(g)\ar[r] & \mathbf{M}\ar[r]^{g} &\mathbf{N}\ar[r]&0
          }$$
commutes.
\end{lem}

\begin{proof}
Without loss of generality, we may assume
$$\mathbf{D}= \quad \cdots\to D^{-2}\to D^{-1}\to D^0\to0\to\cdots.$$
To complete this proof, we need to construct a morphism $\beta=(\beta^k)_{k\in\mathbb{Z}}:\mathbf{D} \to \mathbf{M}$
such that the diagram
$$
\xymatrix@C=0.55cm{
             & \cdots \ar[r]  & D^{-2} \ar@{.>}[dd]\ar@{}[d]_{\beta^{-2}} \ar[rr]^{\delta_\mathbf{D}^{-2}}\ar[dr]^{\alpha^{-2}} & & D^{-1}\ar@{.>}[dd] \ar@{}[d]_{\beta^{-1}}\ar[rr]^{\delta_\mathbf{D}^{-1}}\ar[dr]^{\alpha^{-1}}  && D^0 \ar@{.>}[dd]\ar@{}[d]_{\beta^{0}}\ar[rr] \ar[dr]^{\alpha^{0}} && 0 \ar@{.>}[dd]\ar[dr]\ar[r] & \cdots\\
             & & \cdots \ar[r]  & N^{-2} \ar@{-}[r] & \ar[r]  &  N^{-1}\ar@{-}[r] & \ar[r]       & N^0 \ar@{-}[r] & \ar[r]& N^{1} \ar[r]           & \cdots \\
&\cdots \ar[r] & M^{-2} \ar[rr]_{\delta_\mathbf{M}^{-2}}\ar[ru]_{g^{-2}} && M^{-1} \ar[rr]_{\delta_\mathbf{M}^{-1}}  \ar[ru]_{g^{-1}}                 && M^{0} \ar[rr]_{\delta_\mathbf{M}^{0}}\ar[ru]_{g^{0}} &&  M^{1} \ar[r]\ar[ru]_{g^{1}}      & \cdots
}
$$
commutes.

Let $\beta^k=0$ for $k\geqslant1$. Since $(\ddag)$ is exact,
$0\to \Ker (g^i)\to M^i\overset{g^i }{\longrightarrow} N^i\to0$
is Hom$_\mathcal {A}(X,-)$-exact for each $i\in\mathbb{Z}$ and each object $X\in\mathcal {X}$.
Hence, there exists $\gamma^i:D^i\to M^{i}$ such that $g^{i}\gamma^i=\alpha^i$ for each $i\in \mathbb{Z}$.
Note that
$$\Ker(g)=\,\,\, \cdots\to \Ker (g^{-1})\overset{\delta_\mathbf{M}^{-1} }{\longrightarrow}\Ker (g^{0})\overset{\delta_\mathbf{M}^{0}}{\longrightarrow}\Ker (g^1)\overset{\delta_\mathbf{M}^{1} }{\longrightarrow}\Ker (g^2)\to\cdots$$
is $\mathcal {X}$-acyclic.
As $D^0\in\mathcal {X} $, we have an exact sequence
$$\cdots\to\Hom_\mathcal {A}(D^0,\Ker (g^0))\to\Hom_\mathcal {A}(D^0,\Ker (g^1))\to\Hom_\mathcal {A}(D^0,\Ker (g^2))\to\cdots.$$
Since $g^{1}\delta_\mathbf{M}^0\gamma^0 = \delta_\mathbf{N}^0g^0\gamma^0=\delta_\mathbf{N}^0 \alpha^0=0$,
$-\delta_\mathbf{M}^0\gamma^0\in \Hom_\mathcal{A}(D^0,\Ker (g^1))$.
Now $\delta_\mathbf{M}^1\delta_\mathbf{M}^0\gamma^0=0$,
and thus there exists $\mu^0:D^0\to \Ker (g^0)$ such that $\delta_\mathbf{M}^0\mu^0=-\delta_\mathbf{M}^0\gamma^0$.
Let $\beta^0=\mu^0+\gamma^0$.
Then $\delta_\mathbf{M}^0\beta^0=\delta_\mathbf{M}^0\mu^0+\delta_\mathbf{M}^0\gamma^0=0$
and $g^0\beta^0=g^0\mu^0+g^0\gamma^0=\alpha^0$.

Note that
$$\cdots\to\Hom_\mathcal {A}(D^{-1},\Ker (g^{-1}))\to\Hom_\mathcal {A}(D^{-1},\Ker (g^0))\to\Hom_\mathcal {A}(D^{-1},\Ker (g^1))\to\cdots$$
is also exact.
Since $g^0(\beta^0\delta_\mathbf{D}^{-1}-\delta_\mathbf{M}^{-1}\gamma^{-1})= \alpha^0\delta_\mathbf{D}^{-1}-\delta_\mathbf{N}^{-1}g^{-1}\gamma^{-1}=
\alpha^0\delta_\mathbf{D}^{-1}-\delta_\mathbf{N}^{-1}\alpha^{-1}=0$,
$\beta^0\delta_\mathbf{D}^{-1}-\delta_\mathbf{M}^{-1}\gamma^{-1}\in \Hom_\mathcal {A}(D^{-1},\Ker (g^0))$.
Moreover, since $\delta_\mathbf{M}^0(\beta^0\delta_\mathbf{D}^{-1}-\delta_\mathbf{M}^{-1}\gamma^{-1})
= \delta_\mathbf{M}^0\beta^0\delta_\mathbf{D}^{-1}-\delta_\mathbf{M}^0\delta_\mathbf{M}^{-1}\gamma^{-1}
=\delta_\mathbf{M}^0\beta^0\delta_\mathbf{D}^{-1}=0$,
there exists $\mu^{-1}:D^{-1}\to \Ker (g^{-1})$ such that
$\delta_\mathbf{M}^{-1}\mu^{-1}=\beta^0\delta_\mathbf{D}^{-1}-\delta_\mathbf{M}^{-1}\gamma^{-1}$.
Now let $\beta^{-1}=\mu^{-1}+\gamma^{-1}$, and then
$\delta_\mathbf{M}^{-1}\beta^{-1}=\delta_\mathbf{M}^{-1}\mu^{-1}+\delta_\mathbf{M}^{-1}\gamma^{-1}=\beta^0\delta_\mathbf{D}^{-1}$
and $g^{-1}\beta^{-1}=g^{-1}\mu^{-1}+g^{-1}\gamma^{-1}=\alpha^{-1}$.

Continuing in this manner, we can get $\beta^k=\mu^k+\gamma^k$ such that
$\beta^{k+1}\delta_\mathbf{D}^k=\delta_\mathbf{M}^k\beta^k$ and $g^k\beta^k=\alpha^k$
for $k =-2,-3,\cdots $.
Thus, we obtain a morphism $\beta=(\beta^k)_{k\in\mathbb{Z}}:\mathbf{D}\to\mathbf{M}$
such that $g\beta=\alpha$.
This completes the proof.
\end{proof}

Assume that $\mathcal {X}$ is closed under extensions and has an injective cogenerator.
Then according to \cite[Lemma 3.3(b)]{16}, one has $\res\,\widehat{\mathcal {X}}\subseteq\res\,\widetilde{\mathcal {X}}$, that is,
every object in $\mathcal {A}$ with finite $\mathcal {X}$-projective dimension admits a proper $\mathcal {X}$-resolution.
This fact helps us to get the following result, which plays a key role in the proof of Theorem \ref{thm 3.10}.

\begin{lem}\label{lem 3.9}
Assume that $\mathcal {X}$ is closed under extensions and has an
injective cogenerator. Then for each object $\mathbf{T}\in
\mathrm{\mathbf{K}}^-(\res\, \widehat{\mathcal {X}})$, there
exists a short exact sequence
$$0\to \mathbf{K}\to \mathbf{D}\to \mathbf{T}\to0$$
of complexes such that $\mathbf{D}\in
\mathrm{\mathbf{K}}^{-}(\mathcal {X})$, $\mathbf{K}$ is
$\mathcal{X}$-acyclic and it remains exact after applying the
functor $\Hom_\mathcal {A}(\mathcal {X},-)$, i.e., the sequence
$$0\to\Hom_\mathcal {A}(X,\mathbf{K})\to \Hom_\mathcal {A}(X,\mathbf{D})\to
\Hom_\mathcal {A}(X,\mathbf{T})\to 0 $$ of complexes is
exact for each object $X\in\mathcal {X}$.
\end{lem}

\begin{proof}
Without loss of generality,
we may assume
$$\mathbf{T}= \cdots \to T^{-2} \to T^{-1} \to T^0 \to 0 \to \cdots.$$
Let
$\mathbf{T}(n)= \cdots\to0\to T^{-n}\to T^{-n+1}\to\cdots\to T^0\to 0\to\cdots$
for $n\geqslant 0$.
We will divide the proof into three steps. \vskip 2mm

\textbf{Step 1}. There exists a short exact sequence
$0\to \mathbf{K}(0)\to \mathbf{D}(0)\overset{\underline{\alpha}_0}{\longrightarrow}
               \mathbf{T}(0)\to 0 $
of complexes such that
$\mathbf{D}(0)\in \mathrm{\mathbf{K}}^{-}(\mathcal {X})$,
$\mathbf{K}(0)$ is $\mathcal{X}$-acyclic
and it remains exact after applying the functor $\Hom_\mathcal {A}(X,-)$ for each object $X\in\mathcal {X}$.
\vskip 2mm

Indeed, since $T^0\in\res\, \widehat{\mathcal {X}}$
it admits a proper $\mathcal {X}$-resolution
$$\mathbf{D}(0)=\quad\cdots\to D_0^{-n}\to\cdots\to D_0^{-1}\to D_0^0\to 0\to\cdots$$
by \cite[Lemma 3.3(b)]{16},
i.e., there is an associated exact sequence
$$\mathbf{D}(0)^+=\quad\cdots\to D_0^{-n}\to\cdots\to D_0^{-1}\to D_0^0\overset{\alpha_0}{\longrightarrow} T^0\to0.$$
Hence, we obtain the following short exact sequence of
complexes
$$
\xymatrix{
              & 0        \ar[d]            &               & 0\ar[d]                    & 0              \ar[d]                  & 0 \ar[d]       &        & 0             \ar[d]                        &          \\
\cdots \ar[r] & D_0^{-n} \ar[r]\ar[d]^{\|} & \cdots \ar[r] & D_0^{-1} \ar[r]\ar[d]^{\|} & \Ker(\alpha_0) \ar[r]\ar[d]            & 0 \ar[d]\ar[r] & \cdots & \mathbf{K}(0) \ar[d]                        &          \\
\cdots \ar[r] & D_0^{-n} \ar[r]\ar[d]      & \cdots \ar[r] & D_0^{-1} \ar[r]\ar[d]      & D_0^0          \ar[r]\ar[d]^{\alpha_0} & 0 \ar[r]\ar[d] & \cdots & \mathbf{D}(0) \ar[d]^{\underline{\alpha}_0} & (\sharp) \\
\cdots \ar[r] & 0        \ar[r]\ar[d]      & \cdots \ar[r] & 0        \ar[r]\ar[d]      & T^0            \ar[r]\ar[d]            & 0 \ar[r]\ar[d] & \cdots & \mathbf{T}(0) \ar[d]                        &          \\
              & 0                          &               & 0                          & 0                                      & 0              &        & 0                                           &           }
$$
It is obvious that
$$(\sharp) \quad 0 \to \mathbf{K}(0) \to \mathbf{D}(0) \xrightarrow{\ \underline{\alpha}_0 \ } \mathbf{T}(0) \to 0$$
remains exact after applying the
functor $\Hom_\mathcal {A}(X,-)$ for each object $X\in \mathcal
{X}$ and $\mathbf{K}(0)$ is $\mathcal {X}$-acyclic.
\vskip 2mm

\textbf{Step 2}. For every $n \geqslant 1$, there exists a short exact sequence
$$0\to \mathbf{K}(n) \xrightarrow{\lambda_n} \mathbf{D}(n)\xrightarrow{\ \underline{\alpha}_n \ } \mathbf{T}(n)\to 0 \quad (\natural_n)$$
of complexes such that $\mathbf{D}(n)\in
\mathrm{\mathbf{K}}^{-}(\mathcal {X})$, $\mathbf{K}(n)$ is
$\mathcal{X}$-acyclic and it remains exact after applying the
functor $\Hom_\mathcal {A}(X,-)$ for each object $X\in\mathcal
{X}$.
\vskip 2mm

To see this, we may assume that are done for $n-1$.
By a similar argument as in Step 1,
one has a short exact sequence
$$0 \to \mathbf{L} \to \mathbf{X} \xrightarrow{\ \underline{\beta}\ } T^{-n} \to 0$$
of complexes, where
$\mathbf{L} = \cdots \to X^{-n} \to \cdots \to X^{-1} \to \Ker(\beta) \to 0 \to \cdots$ is $\mathcal {X}$-acyclic
and $\mathbf{X} = \cdots \to X^{-n} \to \cdots \to X^{-1} \to X^0 \to 0 \to \cdots$.
Define a morphism $\sigma: T^{-n}[n-1] \to \mathbf{T}(n-1)$ of complexes as follows:
$$\xymatrix{
\cdots \ar[r] & 0 \ar[r]\ar[d] & T^{-n}   \ar[r]\ar[d]^{\delta_{\mathbf{T}}^{-n}} & 0         \ar[r]\ar[d] & \cdots \ar[r] & 0   \ar[r]\ar[d] & 0 \ar[r]\ar[d] & \cdots \\
\cdots \ar[r] & 0 \ar[r]       & T^{-n+1} \ar[r]^{\delta_{\mathbf{T}}^{-n+1}}     & T^{-n+2}  \ar[r]       & \cdots \ar[r] & T^0 \ar[r]       & 0 \ar[r]       & \cdots }$$
Then it follows from Lemma \ref{lem 3.8} that there exists a morphism
$\psi: \mathbf{X}[n-1] \to \mathbf{D}(n-1)$ of complexes such that the following diagram
$$\xymatrix{
0 \ar[r] & \mathbf{L}[n-1] \ar[r]\ar@.[d] & \mathbf{X}[n-1] \ar[r]\ar@.[d]^{\psi} & T^{-n}[n-1]     \ar[r]\ar[d]^{\sigma} & 0 \\
0 \ar[r] & \mathbf{K}(n-1) \ar[r]         & \mathbf{D}(n-1) \ar[r]                & \mathbf{T}(n-1) \ar[r]                & 0 }$$
commutes.
In particular, the diagram
$$\xymatrix{
 X^{0}          \ar[rr]^{\beta} \ar[d]_{\psi^{-n+1}} & & T^{-n} \ar[d]^{\delta_{\mathbf{T}}^{-n}} \\
 D^{-n+1}_{n-1} \ar[rr]^{\alpha_{n-1}^{-n+1}}        & & T^{-n+1}                                 }$$
commutes.

Let $\mathbf{D}(n)$ be the mapping cone of $\psi: \mathbf{X}[n-1]\to \mathbf{D}(n-1)$, i.e.,
$$\mathbf{D}(n) = \cdots \to                                           X^{-1}\oplus D_{n-1}^{-n-1}
                         \xrightarrow{\delta_{\mathbf{D}(n)}^{-n-1}}   X^0   \oplus D_{n-1}^{-n}
                         \xrightarrow{\delta_{\mathbf{D}(n)}^{-n}}     D^{-n+1}_{n-1}
                         \xrightarrow{\delta_{\mathbf{D}(n-1)}^{-n+1}} D^{-n+2}_{n-1}
                         \to                                           \cdots $$
in which $\delta_{\mathbf{D}(n)}^{-n}     = \left(\begin{array}{cc} 0 & 0 \\ \psi^{-n+1}& \delta_{\mathbf{D}(n-1)}^{-n} \end{array}\right)$
and      $\delta_{\mathbf{D}(n)}^{-n-i+1} = \left(\begin{array}{cc} (-1)^n\delta_{\mathbf{X}}^{-i+2} & 0 \\ \psi^{-n-i+2} & \delta_{\mathbf{D}(n-1)}^{-n-i+1} \end{array}\right)$
for $i \geqslant 2$.
It is trivial that $\mathbf{D}(n)\in \mathrm{\mathbf{K}}^{-}(\mathcal {X})$
because $\mathcal {X}$ is closed under extensions.

We now define a surjective morphism $\mathbf{D}(n) \to \mathbf{T}(n)$ of complexes as follows:
$$\xymatrix{
\cdots \ar[r] & X^{-1}\oplus D_{n-1}^{-n-1} \ar[r]^{\delta_{\mathbf{D}(n)}^{-n-1}}\ar[d]^{0} & X^0\oplus D_{n-1}^{-n} \ar[r]^{\quad \delta_{\mathbf{D}(n)}^{-n}} \ar[d]^{(\beta, 0)} & D^{-n+1}_{n-1} \ar[r]^{\delta_{\mathbf{D}(n-1)}^{-n+1}} \ar[d]^{\alpha_{n-1}^{-n+1}} & D^{-n+2}_{n-1} \ar[r]\ar[d]^{\alpha_{n-1}^{-n+2}} & \cdots \\
\cdots \ar[r] & 0                           \ar[r]                                           & T^{-n}                 \ar[r]^{\delta_{\mathbf{T}}^{-n}}                              & T^{-n+1}       \ar[r]^{\delta_{\mathbf{T}}^{-n+1}}                                   & T^{-n+2}       \ar[r]                             & \cdots }$$
Then we obtain a short exact sequence
$$0 \to \mathbf{K}(n) \to \mathbf{D}(n) \to \mathbf{T}(n) \to 0,$$
of complexes, which remains exact after applying the functor \Hom$_\mathcal {A}(X,-)$ for each $X\in \mathcal{X}$.

Define a morphism $\varphi: \mathbf{K}(n-1)\to\mathbf{K}(n)$ of complexes as follows:
$$\xymatrix{
\cdots \ar[r] & D_{n-1}^{-n-1}               \ar[r]\ar[d]^{(0,1)} & D_{n-1}^{-n}                    \ar[r]\ar[d]^{(0,1)} & \Ker(\alpha^{-n+1}_{n-1}) \ar[r]\ar[d]^{\|} & \cdots \\
\cdots \ar[r] & X^{-1} \oplus D_{n-1}^{-n-1} \ar[r]               & \Ker(\beta) \oplus D_{n-1}^{-n} \ar[r]               & \Ker(\alpha^{-n+1}_{n-1}) \ar[r]            & \cdots }$$
Then we get a short exact sequence
$$0 \to \mathbf{K}(n-1) \xrightarrow{\ \varphi\ } \mathbf{K}(n) \to \mathbf{L}[n] \to 0$$
of complexes,
which remains exact after applying the functor $\Hom_\mathcal {A}(X,-)$ for each object $X\in \mathcal {X}$.
Consequently, $\mathbf{K}(n)$ is also $\mathcal{X}$-acyclic
since $\mathbf{K}(n-1)$ is so by the induction hypothesis.

\vskip 2mm
\textbf{Step 3}. Set
$$\mathbf{D}(n) = \cdots \to D_n^{-n} \xrightarrow{\delta_{\mathbf{D}(n)}^{-n}}
                             D_n^{-n+1} \to \cdots \to
                             D_n^{-1} \xrightarrow{\delta_{\mathbf{D}(n)}^{-1}}
                             D_n^0 \to 0 \to \cdots$$
and
$$\mathbf{K}(n) = \cdots \to K_n^{-n} \xrightarrow{\delta_{\mathbf{K}(n)}^{-n}}
                             K_n^{-n+1} \to \cdots \to
                             K_n^{-1} \xrightarrow{\delta_{\mathbf{K}(n)}^{-1}}
                             K_n^0 \to 0 \to \cdots$$
for each $n \geqslant 0$.
Note that $D_n^{-n+1} = D_{n-1}^{-n+1}$ and $K_n^{-n+1} = K_{n-1}^{-n+1}$ for all $n \geqslant 1$.
So we have the following complexes
$$\mathbf{D} = \cdots \to D_n^{-n} \xrightarrow{\delta_{\mathbf{D}}^{-n}}
                          D_{n-1}^{-n+1} \to \cdots \to
                          D_1^{-1} \xrightarrow{\delta_{\mathbf{D}}^{-1}}
                          D_0^0 \to 0 \to \cdots$$
and
$$\mathbf{K} = \cdots \to K_n^{-n} \xrightarrow{\delta_{\mathbf{K}}^{-n}}
                          K_{n-1}^{-n+1} \to \cdots \to
                          K_1^{-1} \xrightarrow{\delta_{\mathbf{K}}^{-1}}
                          K_0^0 \to 0 \to \cdots$$
where $\delta_{\mathbf{D}}^{-n} = \delta_{\mathbf{D}(n)}^{-n}$
and $\delta_{\mathbf{K}}^{-n} = \delta_{\mathbf{K}(n)}^{-n}$ for each $n \geqslant 0$.
In addition, there exists a short exact sequence
$$0 \to \mathbf{K} \xrightarrow{\lambda} \mathbf{D} \xrightarrow{\alpha} \mathbf{T}\to 0$$
of complexes,
where $\lambda^{-n} = \lambda^{-n}_n: K_n^{-n} \to D_n^{-n}$ and
$\alpha^{-n} = \underline{\alpha}_n^{-n}: D_n^{-n} \to T_n^{-n}$
for each $n \geqslant 0$ (see Step 2).

It is easy to see that
$0 \to \mathbf{K} \xrightarrow{\lambda} \mathbf{D} \xrightarrow{\alpha} \mathbf{T}\to 0$
is as desired.
This completes the proof.
\end{proof}

From \cite[Theorem 4.6(1)]{Huang} we know that
$\res\,\widetilde{\mathcal {X}}$ is closed under direct summands.
This fact enables us to get the next result.

\begin{lem}
Assume that $\mathcal {X}$ is exact and has an injective cogenerator.
Then $\res\,\widehat{\mathcal {X}}$ is closed under direct summands.
\end{lem}
\begin{proof}
Let $M$ be an object in $\mathcal {A}$ such that $M\in\res\,\widehat{\mathcal {X}}$
and $M'$ any direct summand of $M$.
It suffices to prove that $M'\in\res\,\widehat{\mathcal {X}}$.
Note that $M$ is also in $\res\,\widetilde{\mathcal {X}}$ (see \cite[Lemma 3.3(b)]{16}).
It follows from \cite[Theorem 4.6(1)]{Huang} that $M'\in \res\,\widetilde{\mathcal {X}}$.
Since $\mathcal {X}$ is closed under extensions and has an injective cogenerator,
Ext$_{\mathcal {XA}}^n(M,-)=0$ for $n>\mathcal {X}$-pd$(M)$ by \cite[Proposition 4.5(b)(2)]{16}.
Hence, Ext$_{\mathcal {XA}}^n(M',-)=0$ for $n>\mathcal {X}$-pd$(M)$,
and so $M'\in\res\,\widehat{\mathcal {X}}$ by \cite[Proposition 4.5(a)]{16}
since $\mathcal {X}$ is also closed under direct summands, as desired.
\end{proof}

Now let $H: \mathrm{\mathbf{K}}^{-}(\mathcal {X})\to\mathrm{\mathbf{D}}^{-}_{\mathcal{X}}(\res\,\widehat{\mathcal {X}})$
be the composition of the embedding $\mathrm{\mathbf{K}}^{-}(\mathcal {X})\hookrightarrow\mathrm{\mathbf{K}}^{-}(\res\,\widehat{\mathcal {X}})$
and the localization
$Q:\mathrm{\mathbf{K}}^{-}(\res\,\widehat{\mathcal {X}})\to\mathrm{\mathbf{D}}^{-}_{\mathcal{X}}(\res\,\widehat{\mathcal {X}})$.
Clearly, $H$ is a triangle functor.
Moreover, Lemma \ref{lem 3.7} implies that $H$ is fully faithful (the case $\mathcal {S}=\res\,\widehat{\mathcal {X}}$),
and from Lemma \ref{lem 3.9} we know that $H$ is dense.
Thus, we obtain the main result of this paper.

\begin{thm}\label{thm 3.10}
Assume that $\mathcal {X}$ is exact and has an injective cogenerator.
Then there exists a triangle-equivalence
$\mathrm{\mathbf{K}}^{-}(\mathcal {X})\cong \mathrm{\mathbf{D}}^{-}_{\mathcal{X}}(\res\,\widehat{\mathcal {X}})$.
\end{thm}

In \cite[Section 3]{4} Chen investigated the relative derived category of the case $\mathcal {S}=\mathcal {A}$ in Definition \ref{def 3.4}.
Let $\mathcal {X}\subseteq\mathcal {A}$ be a contravariantly finite subcategory.
Assume that $\mathcal {X}$ is admissible and $\mathcal {X}$-pd$(\mathcal {A}) < \infty$.
Chen proved in \cite[Propositin 3.5]{4} that the natural composite functor
$\mathrm{\mathbf{K}}(\mathcal {X})\hookrightarrow\mathrm{\mathbf{K}}(\mathcal{A})
\overset{Q}{\longrightarrow}\mathrm{\mathbf{D}}_{\mathcal{X}}(\mathcal {A})$
is a triangle-equivalence.

In view of the following example, we can conclude that Theorem \ref{thm 3.10} can not be deduced from \cite[Propositin 3.5]{4}.

\begin{exa}\label{exa 3.10}
Let $R=\mathbb{Z}_4$ and $\mathcal {X}=\mathcal {P}_R$,
where $\mathcal {P}_R$ denotes the class of all projective $R$-modules.
It is clear that $\mathcal {X}$ is exact
with itself as an injective cogenerator.
Then it follows from {\rm Theorem \ref{thm 3.10}} that there exists a triangle-equivalence
$\mathrm{\mathbf{K}}^{-}(\mathcal {X})\cong \mathrm{\mathbf{D}}^{-}_{\mathcal{X}}(\res\,\widehat{\mathcal {X}})$.
But here $\res\,\widehat{\mathcal {X}}$ is not an abelian category
since the kernel of the morphism
$f:\mathbb{Z}_4\to \mathbb{Z}_4$ with $f(x)=2\,x$ is $2\,\mathbb{Z}_4$ and
$2\,\mathbb{Z}_4$ has infinite projective dimension.
\end{exa}

Recall from \cite[Definition 4.1]{17} that an acyclic complex of objects in $\mathcal {X}$ is called $totally\,\mathcal {X}$-$acyclic$
if it is $\Hom_\mathcal {A}(\mathcal {X},-)$-exact and $\Hom_\mathcal {A}(-,\mathcal {X})$-exact.
Let $\mathcal {G}(\mathcal {X})$ denote the subcategory of $\mathcal {A}$ whose modules are of the form $M \cong\Ker(\delta^{1}_\mathbf{X})$
for some totally $\mathcal {X}$-acyclic complex $\mathbf{X}$.
In the special case $\mathcal {X}=\mathcal {P}$,
the objects in $\mathcal {G}(\mathcal {P})$
are called Gorenstein projective objects of $\mathcal {A}$.

\begin{rem}\label{rem 3.12}
Assume $\mathcal {X}\perp\mathcal {X}$.
Then by virtue of \cite[Theorem B]{17},
$\mathcal {G}(\mathcal {X})$ is an exact subcategory of $\mathcal {A}$, and
it is closed under kernels of epimorphisms if $\mathcal {X}$ is so.
Moreover, it follows from \cite[Corollary 4.7]{17} that
$\mathcal {X}$ is both an injective cogenerator and a projective generator for $\mathcal {G}(\mathcal {X})$,

\end{rem}

According to Remark \ref{rem 3.12}, we can obtain the next result by Theorem \ref{thm 3.10}.

\begin{cor}\label{cor 3.13}
Assume $\mathcal {X}\perp\mathcal {X}$.
Then there exists a triangle-equivalence
$\mathrm{\mathbf{K}}^{-}(\mathcal {G(X)})\\ \cong \mathrm{\mathbf{D}}^{-}_{\mathcal{G(X)}}(\res\,\widehat{\mathcal {G(X)}})$.
\end{cor}

It is obvious that $\mathcal {X}=\mathcal {P}$ satisfies the hypothesis of Corollary \ref{cor 3.13}.
Hence, we have the next triangle-equivalence, which extends \cite[Theorem 3.6(1)]{7} to the bounded below case.

\begin{cor}\label{cor 3.14}
$\mathrm{\mathbf{K}}^{-}(\mathcal {G(P)})\cong \mathrm{\mathbf{D}}^{-}_{\mathcal{G(P)}}(\res\,\widehat{\mathcal {G(P)}})$.
\end{cor}

The notion of semidualizing module goes back to \cite{5} (where the more general PG-module are studied) and \cite{6}.
Christensen \cite{2} extended this notion to semidualizing complex.
Recently, it has already been defined over arbitrary associative rings \cite[Definition 2.1]{9}.
In what follows we assume that $R$ is a commutative ring.
Recall from \cite{19} that an $R$-module $C$ is called $semidualizing$ if
$C$ admits a degree-wise finite projective resolution, $\Ext^{\geqslant1}_R (C,C) = 0$,
and the natural homothety map $R\to\Hom_R(C,C)$ is an isomorphism.
More examples of semidualizing module can be found in \cite{2, 9, 10}.

Following \cite{19}, an $R$-module is called $C$-$projective$ (resp., $C$-$flat$)
if it is isomorphic to an $R$-module of the form $C\otimes_R P$ for some projective (resp., flat) $R$-module $P$.
An $R$-module is called $C$-\emph{injective} if it is isomorphic to an $R$-module of the form Hom$_R(C,I)$ for some injective $R$-module $I$.
Let $\mathcal {P}_C$ (resp., $\mathcal {F}_C$, $\mathcal {I}_C$)
denote the subcategory of $C$-projective (resp., $C$-flat, $C$-injective) $R$-modules.
A $complete$ $\mathcal {P}\mathcal {P}_C$-$resolution$ is an exact and $\Hom_R(-,\mathcal {P}_C)$-exact complex $\mathbf{X}$
of $R$-modules with $X_i$ projective for $i \leqslant 0$ and $X_i$ $C$-projective for $i > 0$.
An $R$-module $M$ is $\mathcal {G}_C$-$projective$ if
there exists a complete $\mathcal {P}\mathcal {P}_C$-resolution $\mathbf{X}$ such that $M \cong\Ker(\delta^{1}_\mathbf{X})$.
A $complete$ $\mathcal {F}\mathcal {F}_C$-$resolution$ is an exact and $-\otimes_R \mathcal {I}_C$-exact complex $\mathbf{Z}$
of $R$-modules with $Z_i$ flat for $i \leqslant 0$ and $Z_i$ $C$-flat for $i > 0$.
An $R$-module $T$ is $\mathcal {G}_C$-$flat$ if
there exists a complete $\mathcal {F}\mathcal {F}_C$-resolution $\mathbf{Z}$ such that $T \cong\Ker(\delta^{1}_\mathbf{Z})$.
Let $\mathcal {GP}_C$ (resp., $\mathcal {GF}_C$) denote
the subcategory of $\mathcal {G}_C$-projective (resp., $\mathcal {G}_C$-flat) $R$-modules.

For more details about semidualizing modules and their related categories,
we refer the reader to [12, 14, 19-24].

\begin{cor}\label{cor 3.14}
Let $R$ be a commutative ring and $C$ a semidualizing $R$-module.
Then there exists a triangle-equivalence
$\mathrm{\mathbf{K}}^{-}(\mathcal {GP}_C)\cong \mathrm{\mathbf{D}}^{-}_{\mathcal {GP}_C}(\res\,\widehat{\mathcal {GP}_C})$.
\end{cor}

\begin{proof}
From \cite[Proposition 3.15]{19} we know that $\mathcal {G}\mathcal{P}_C$ is exact.
Moreover, according to \cite[Propsition 2.6 and 2.7]{19},
$\mathcal {P}_C\subseteq\mathcal {G}\mathcal{P}_C$ and $\mathcal {G}\mathcal {P}_C\perp\mathcal {P}_C$.
Hence, it follows from \cite[Proposition 3.16]{19} that $\mathcal {P}_C$ is an injective cogenerator for $\mathcal {G}\mathcal{P}_C$.
Then the result follows by Theorem \ref{thm 3.10}.
\end{proof}

If $R$ is commutative noetherian,
\cite[Theorem $\textrm{I}$]{15} said that $\mathcal {GF}_C$ is exact
and $\mathcal {F}^{\textrm{cot}}_C$ is an injective cogenerator for $\mathcal {GF}_C$,
where $\mathcal {F}^{\textrm{cot}}_C$ denotes the subcategory of $R$-modules $C\otimes_R H$ with $H$ flat and cotorsion.
Then by Theorem \ref{thm 3.10}, we have

\begin{cor}\label{cor 3.17}
Let $R$ be a commutative noetherian ring and $C$ a semidualizing $R$-module. Then there exists a triangle-equivalence $\mathrm{\mathbf{K}}^{-}(\mathcal {GF}_C)\cong \mathrm{\mathbf{D}}^{-}_{\mathcal {GF}_C}(\res\,\widehat{\mathcal {GF}_C})$.
\end{cor}

\section{\bf Relative Cohomology and Applications}

We begin this section with the following result, which will be used in the proofs of Theorem \ref{thm 4.2} and Proposition \ref{prop 4.6}.

\begin{prop}\label{prop 4.1}
$\mathrm{\mathbf{D}}^{-}_{\mathcal{X}}(\mathcal {A})$ is a full triangulated subcategory of $\mathrm{\mathbf{D}}_{\mathcal{X}}(\mathcal {A})$.
$\mathrm{\mathbf{D}}^{b}_{\mathcal{X}}(\mathcal {A})$ is a full triangulated subcategory of $\mathrm{\mathbf{D}}^{-}_{\mathcal{X}}(\mathcal {A})$, and
hence of $\mathrm{\mathbf{D}}_{\mathcal{X}}(\mathcal {A})$.
\end{prop}

\begin{proof}
We only prove the first assertion, and the second one can be proved similarly.
Let $\Sigma_\mathcal {X}$ be the compatible multiplicative system determined by $\mathrm{\mathbf{K}}_{\mathcal{X}}(\mathcal {A})$, that is,
the collection of all $\mathcal{X}$-quasi-isomorphisms in $\mathrm{\mathbf{K}}(\mathcal {A})$. Then
$\mathrm{\mathbf{D}}_{\mathcal{X}}(\mathcal {A}) = \Sigma_\mathcal {X}^{-1}\mathrm{\mathbf{K}}(\mathcal {A})$ and
$\mathrm{\mathbf{D}}^{-}_{\mathcal{X}}(\mathcal {A}) = (\Sigma_\mathcal {X}\cap \mathrm{\mathbf{K}}^{-}(\mathcal {A}))^{-1}\, \mathrm{\mathbf{K}}^{-}(\mathcal {A})$.
By \cite[Proposition (III) 2.10]{8} it suffices to prove that for any $\mathcal{X}$-quasi-isomorphism
$f: \mathbf{X}\to \mathbf{Y}$ with $\mathbf{Y} \in \mathrm{\mathbf{K}}^{-}(\mathcal {A})$,
there is a morphism $g: \mathbf{X}^{'}\to \mathbf{X}$ with $\mathbf{X}^{'}\in
\mathrm{\mathbf{K}}^{-}(\mathcal {A})$ such that $fg$ is an $\mathcal{X}$-quasi-isomorphism.
Then the canonical functor
$(\Sigma_\mathcal {X}\cap \mathrm{\mathbf{K}}^{-}(\mathcal {A}))^{-1}\, \mathrm{\mathbf{K}}^{-}(\mathcal {A})\to \Sigma_\mathcal {X}^{-1}\mathrm{\mathbf{K}}(\mathcal {A})$
is fully faithful, and hence
$\mathrm{\mathbf{D}}^{-}_{\mathcal{X}}(\mathcal {A})$ is a full triangulated subcategory of $\mathrm{\mathbf{D}}_{\mathcal{X}}(\mathcal {A})$.

Suppose that there is an integer $i$ such that $Y^{k} = 0$ for any $k > i$.
Let $\mathbf{X}^{'}$ be the soft truncation $\mathbf{X}_{{i}{\supset}}$ of $\mathbf{X}$.
Then there is a commutative diagram
\begin{center}
$\xymatrix@C=25pt{
\mathbf{X}_{{i}{\supset}} \ar[d]_{g} &\cdots \ar[r] &X^{i-2}\ar[r]\ar@{=}[d] &X^{i-1}\ar[r]\ar@{=}[d] &
\mathrm{Ker}(\delta_\mathbf{X}^{i})\ar[r]\ar@{^{(}->}[d]&0\ar[r]\ar[d]&\cdots\\
\mathbf{X} \ar[d]_{f}&\cdots \ar[r] &X^{i-2}\ar[r]\ar[d] &X^{i-1}\ar[r]\ar[d] & X^{i}\ar[r]\ar[d]&X^{i+1}\ar[r]\ar[d]&\cdots\\
\mathbf{Y} &\cdots \ar[r] &Y^{i-2}\ar[r]&Y^{i-1}\ar[r]&Y^{i}\ar[r]&0\ar[r]&\cdots\\
}$
\end{center}
It is easy to check that $g$ is also an $\mathcal{X}$-quasi-isomorphism,
and so is $fg$. This completes the proof.
\end{proof}

Now we are in a position to give our anther main result, which is an analog of \cite[Theorem 3.12]{7}.

\begin{thm}\label{thm 4.2}
Let $M$ and $N$ be objects in $\mathcal {A}$.
Assume that $M\in\res\,\widetilde{\mathcal {X}}$.
Then $\mathrm{Ext}^{n}_{\mathcal{X}\mathcal {A}}(M, N) = \Hom_{\mathrm{\mathbf{D}}^{b}_{\mathcal{X}}(\mathcal {A})}(M, N[n])$.
\end{thm}

\begin{proof}
Let $\xymatrix@C=0.5cm{\mathbf{X} \ar[r]^{\simeq} & M}$ be a proper $\mathcal{X}$-resolution of $M$.
Then $\xymatrix@C=0.5cm{ \mathbf{X} \ar[r]^{\simeq} & M}$ is an $\mathcal{X}$-quasi-isomorphism,
and so $\mathbf{X} \cong M$ in $\mathrm{\mathbf{D}}^{-}_{\mathcal{X}}(\mathcal {A})$.
This yields the third isomorphism in the following sequence
$$\begin{aligned}
\mathrm{Ext}^{n}_{\mathcal{X}\mathcal {A}}(M, N) &= \mathrm{H}^{n}\,(\Hom_{\mathcal {A}}(\mathbf{X}, N))\\
&\cong \Hom_{\mathrm{\mathbf{K}}^{-}(\mathcal {A})}(\mathbf{X}, N[n])\\
&\cong \Hom_{\mathrm{\mathbf{D}}_{\mathcal{X}}^{-}(\mathcal {A})}(\mathbf{X}, N[n])\\
&\cong \Hom_{\mathrm{\mathbf{D}}^{-}_{\mathcal{X}}(\mathcal {A})}(M, N[n])\\
&\cong \Hom_{\mathrm{\mathbf{D}}^{b}_{\mathcal{X}}(\mathcal {A})}(M, N[n]),
\end{aligned}$$
while the second isomorphism follows from Lemma \ref{lem 3.7} (the case $\mathcal {S}=\mathcal {A}$),
and the fourth one holds by Proposition \ref{prop 4.1}.
This completes the proof.
\end{proof}

The next result follows directly from \cite[Lemma 3.3(b)]{16} and Theorem \ref{thm 4.2}.
It will be used in the proof of Proposition \ref{prop 4.6}.

\begin{cor}\label{cor 4.4}
Let $M$ be an object in $\mathcal {A}$ such that $M\in \res\,\widehat{\mathcal {X}}$.
Assume that $\mathcal {X}$ is closed under extensions and has an injective cogenerator.
Then $\mathrm{Ext}^{n}_{\mathcal{X}\mathcal {A}}(M, N) =\Hom_{\mathrm{\mathbf{D}}^{b}_{\mathcal{X}}(\mathcal {A})}(M, N[n])$.
\end{cor}

In the remainder of the paper, we display two applications of Theorem \ref{thm 4.2}.
Firstly, we give in the relative derived category a brief proof for
the results of long exact sequences about relative derived functor $\mathrm{Ext}_{\mathcal{X}\mathcal {A}}(-, -)$.
For the similar results with different methods, we refer the reader to \cite[Lemma 4.3]{16} or \cite[Lemma 8.2.3 and 8.2.5]{EJ2000}.

\begin{prop}\label{prop 4.5}
Let $M$ be an object in $\mathcal {A}$ and
$N^{\bullet} = \, 0\to N\stackrel{f}\longrightarrow N^{'}\stackrel{g}\longrightarrow N^{''}\to 0$
an $\mathcal{X}$-acyclic sequence in $\mathcal {A}$.

$(1)$ If $M\in\res\,\widetilde{\mathcal {X}}$,
then there exist homomorphisms $\vartheta^{n}_{\mathcal{X}}(M, N^{\bullet})$,
which are natural in $M$ and $N^{\bullet}$,
such that the sequence below is exact
\begin{center}
$\xymatrix@C=50pt{
\cdots\ar[r] &\mathrm{Ext}_{\mathcal{X}\mathcal{A}}^{n}(M, N)\ar[r]^{\mathrm{Ext}_{\mathcal{X}\mathcal{A}}^{n}(M, f)}  &\mathrm{Ext}_{\mathcal{X}\mathcal{A}}^{n}(M, N^{'}) \ar[r]^{\mathrm{Ext}_{\mathcal{X}\mathcal{A}}^{n}(M, g)} &\mathrm{Ext}_{\mathcal{X}\mathcal{A}}^{n}(M, N^{''})\\
 \ar[r]^{\vartheta^{n}_{\mathcal{X}}(M, N^{\bullet})} &\mathrm{Ext}_{\mathcal{X}\mathcal{A}}^{n+1}(M, N) \ar[r]^{\mathrm{Ext}_{\mathcal{X}\mathcal{A}}^{n+1}(M, f)} &\mathrm{Ext}_{\mathcal{X}\mathcal{A}}^{n+1}(M, N^{'}) \ar[r]^{\mathrm{Ext}_{\mathcal{X}\mathcal{A}}^{n+1}(M, g)} &\cdots
.}$
\end{center}

$(2)$ If $N$, $N'$ and $ N^{''}\in\res\,\widetilde{\mathcal {X}}$,
then there exist homomorphisms $\vartheta^{n}_{\mathcal{X}}(N^{\bullet}, M)$,
which are natural in $N^{\bullet}$ and $M$,
such that the sequence below is exact
\begin{center}
$\xymatrix@C=50pt{
\cdots\ar[r] &\mathrm{Ext}_{\mathcal{X}\mathcal{A}}^{n}(N^{''}, M)\ar[r]^{\mathrm{Ext}_{\mathcal{X}\mathcal{A}}^{n}(g, M)}  &\mathrm{Ext}_{\mathcal{X}\mathcal{A}}^{n}(N^{'}, M) \ar[r]^{\mathrm{Ext}_{\mathcal{X}\mathcal{A}}^{n}(f, M)} &\mathrm{Ext}_{\mathcal{X}\mathcal{A}}^{n}(N, M)\\
 \ar[r]^{\vartheta^{n}_{\mathcal{X}}(N^{\bullet}, M)} &\mathrm{Ext}_{\mathcal{X}\mathcal{A}}^{n+1}(N^{''}, M) \ar[r]^{\mathrm{Ext}_{\mathcal{X}\mathcal{A}}^{n+1}(g, M)} &\mathrm{Ext}_{\mathcal{X}\mathcal{A}}^{n+1}(N^{'}, M) \ar[r]^{\mathrm{Ext}_{\mathcal{X}\mathcal{A}}^{n+1}(f, M)}   &\cdots
.}$
\end{center}
\end{prop}

\begin{proof}
Consider all objects $M,\, N,\, N^{'}$ and $N^{''}$ being complexes concentrated in 0th degree,
and $f$, $g$ being morphisms of such complexes.
Since $N^{\bullet}$ is $\mathcal{X}$-acyclic,
$g$ induces an $\mathcal{X}$-quasi-isomorphism $\mathrm{cone}(f)\to N^{''}$ in $\mathrm{\mathbf{K}}^{b}(\mathcal {A})$,
where the mapping cone $\mathrm{cone}(f) =\, \cdots\to 0\to N\stackrel{f}\to N^{'}\to 0\to \cdots$.
Hence, $\mathrm{cone}(f)\cong N^{''}$ in $\mathrm{\mathbf{D}}^{b}_{\mathcal{X}}(\mathcal {A})$.

Consider the following commutative diagram in $\mathrm{\mathbf{D}}^{b}_{\mathcal{X}}(\mathcal {A})$
\begin{center}
$\xymatrix@C=40pt{
N\ar[r]^{f} \ar@{=}[d] &N^{'}\ar[r]\ar@{=}[d] &\mathrm{cone}(f)\ar[r]\ar[d] &N[1]\ar@{=}[d]\\
N\ar[r]^{f} &N^{'}\ar[r]^{g} &N^{''}\ar[r] &N[1]
.}$
\end{center}
This implies that $N\stackrel{f}\to N^{'}\stackrel{g}\to N^{''}\to N[1]$
is a distinguished triangle in $\mathrm{\mathbf{D}}^{b}_{\mathcal{X}}(\mathcal {A})$.
Applying the functors
$\Hom_{\mathrm{\mathbf{D}}^{b}_{\mathcal{X}}(\mathcal {A})}(M, -)$ and $\Hom_{\mathrm{\mathbf{D}}^{b}_{\mathcal{X}}(\mathcal {A})}(-, M)$
respectively to this distinguished triangle, and by Theorem \ref{thm 4.2}, we get our desired long exact sequences.\end{proof}

In what follows,
we pay our attention to $_{R}\mathcal {M}$ (the category of left $R$-modules),
and then let $\mathcal {X}$ and $\mathcal {W}$ denote subcategories of $_{R}\mathcal {M}$.

Recall from \cite[Definition 3.1]{18} that a $Tate$ $\mathcal {W}$-$resolution$ of an $R$-module $M$
is a diagram $\mathbf{T}\stackrel{\nu}\to \mathbf{W}\stackrel{\pi}\to M$ of morphisms of complexes,
where $\mathbf{W}\stackrel{\pi}\to M$ is a proper $\mathcal {W}$-resolution,
$\mathbf{T}$ is a totally $\mathcal {W}$-acyclic complex,
and $\nu_{i}$ is bijective for all $i \ll 0$.
For any $R$-module $N$ and each $n\in \mathbb{Z}$,
the $n$th $Tate~cohomology$ $group$ is defined as
$$\widehat{\textrm{Ext}}_{\mathcal {W}\mathcal {M}}^{n}(M, N): = \textrm{H}^{n}\,(\textrm{Hom}_R(\mathbf{T}, N)),$$
which is independent of choices of resolutions and lifting (see \cite[Definition 4.1 and Lemma 3.8]{18}).

In \cite{11} Liang and Yang gave us another way to calculate such Tate cohomology group.
Let $M$ be an object in $\res\,\widehat{\mathcal {X}}$. If $\mathcal {X}$ is exact and closed under kernels
of epimorphisms, and $\mathcal {W}$ is both an injective cogenerator and a projective
generator for $\mathcal {X}$ and closed under direct summands,
then \cite[Lemma 3.4]{18} implies that $M$ admits a Tate $\mathcal {W}$-resolution
$\mathbf{T}\stackrel{\nu}\to \mathbf{W}\stackrel{\pi}\to M$.
Moreover, $M$ admits a proper $\mathcal {X}$-resolution $\xymatrix@C=0.5cm{\mathbf{X }\ar[r]^{\simeq} & M}$ by \cite[Lemma 3.3(b)]{16}.
Let $f: \mathbf{W}\to \mathbf{X}$ be a lifting of the identity Id$_M : M\to M$.
Liang and Yang argued that, for $n\geqslant1$, $\widehat{\textrm{Ext}}_{\mathcal {W}\mathcal {M}}^{n}(M, N)$ is exactly
$\textrm{H}^{n+1}\,(\textrm{Hom}_R(\mathrm{cone}(f), N))$ (see \cite[Proposition 3.13]{11}).
Hence, $\widehat{\textrm{Ext}}_{\mathcal {W}\mathcal {M}}^{n}(M, N)\cong\Hom_{\mathrm{\mathbf{K}}^{-}(R)}(\mathrm{cone}(f)[-n-1],N)$ for $n\geqslant1$.

Next, we give another application of Theorem \ref{thm 4.2}.
Note that the following Avramov-Martsinkovsky type exact sequence follows directly from \cite[Theorem 4.10]{18}.
But here we give also a proof in the relative derived category.

\begin{prop}\label{prop 4.6}
Assume that $\mathcal {X}$ is exact and closed under kernels of epimorphisms,
and $\mathcal {W}$ is both an injective cogenerator and a projective generator for $\mathcal {X}$
and closed under direct summands.
Let $M$ and $N$ be $R$-modules with $M\in\res\,\widehat{\mathcal {X}}$.
Set $d=\mathcal {X}$-$\pd (M)$.
Then there exists an exact sequence
$$0\longrightarrow \mathrm{Ext}_{\mathcal{X}\mathcal{M}}^{1}(M,
N)\longrightarrow\mathrm{Ext}_{\mathcal{W}\mathcal{M}}^{1}(M,
N)\longrightarrow\widehat{\mathrm{Ext}}_{\mathcal{W}\mathcal{M}}^{1}(M, N)\indent$$
$$\indent\longrightarrow\mathrm{Ext}^{2}_{\mathcal{X}\mathcal {M}}(M,
N)\longrightarrow\mathrm{Ext}_{\mathcal{W}\mathcal{M}}^{2}(M, N)
\longrightarrow\widehat{\mathrm{Ext}}_{\mathcal{W}\mathcal{M}}^{2}(M, N)\longrightarrow$$
$$\indent\cdots\longrightarrow\mathrm{Ext}^{d}_{\mathcal{X}\mathcal {M}}(M,
N)\longrightarrow\mathrm{Ext}_{\mathcal{W}\mathcal{M}}^{d}(M, N)\longrightarrow\widehat{\mathrm{Ext}}_{\mathcal{W}\mathcal{M}}^{d}(M, N)\longrightarrow0.
$$
\end{prop}

\begin{proof}
For $f: \mathbf{W}\to \mathbf{X}$,
there is a distinguished triangle
$$\mathbf{W} \stackrel{f}\longrightarrow \mathbf{X} \stackrel{}\longrightarrow \mathrm{cone}(f)
\stackrel{}\longrightarrow \mathbf{W}[1]$$ in
$\mathrm{\mathbf{K}}^{-}(R)$. By applying the functor
$\Hom_{\mathrm{\mathbf{K}}^{-}(R)}(-, N)$ to it, we then get an
exact sequence
$$\cdots\longrightarrow\Hom_{\mathrm{\mathbf{K}}^{-}(R)}(\mathbf{W}[-n+1],
N)\longrightarrow\Hom_{\mathrm{\mathbf{K}}^{-}(R)}(\mathrm{cone}(f)[-n],N)$$
$$\longrightarrow\Hom_{\mathrm{\mathbf{K}}^{-}(R)}(\mathbf{X}[-n],N)
\longrightarrow\Hom_{\mathrm{\mathbf{K}}^{-}(R)}(\mathbf{W}[-n],N)$$
$$\longrightarrow\Hom_{\mathrm{\mathbf{K}}^{-}(R)}(\mathrm{cone}(f)[-n-1],N)
\longrightarrow\Hom_{\mathrm{\mathbf{K}}^{-}(R)}(\mathbf{X}[-n-1],
N)\longrightarrow\cdots.$$

As $\xymatrix@C=0.5cm{\mathbf{ X }\ar[r]^{\simeq} & M}$ is an $\mathcal{X}$-quasi-isomorphism,
$\mathbf{X} \cong M$ in $\mathrm{\mathbf{D}}^{-}_{\mathcal{X}}(R)$.
This yields the third isomorphism in the following sequence
$$\begin{aligned}
\Hom_{\mathrm{\mathbf{K}}^{-}(R)}(\mathbf{X}[-n], N)
&\cong\Hom_{\mathrm{\mathbf{K}}^{-}(R)}(\mathbf{X}, N[n])\\
&\cong \Hom_{\mathrm{\mathbf{D}}_{\mathcal{X}}^{-}(R)}(\mathbf{X}, N[n])\\
&\cong \Hom_{\mathrm{\mathbf{D}}^{-}_{\mathcal{X}}(R)}(M, N[n])\\
&\cong \Hom_{\mathrm{\mathbf{D}}^{b}_{\mathcal{X}}(R)}(M, N[n])\\
&\cong \mathrm{Ext}_{\mathcal{X}\mathcal{M}}^{n}(M, N),
\end{aligned}$$
while the second isomorphism follows from Lemma \ref{lem 3.7} (the case $\mathcal {S}=R$),
the fourth isomorphism is from Proposition \ref{prop 4.1},
and the last one holds by Corollary \ref{cor 4.4}.
Similarly, we have $\Hom_{\mathrm{\mathbf{K}}^{-}(R)}(\mathbf{W}[-n], N)\cong
\mathrm{Ext}_{\mathcal {W}\mathcal {M}}^{n}(M, N)$.
By left exactness of $\Hom$, it is easy to check that
$\Hom_{\mathrm{\mathbf{K}}^{-}(R)}(\mathrm{cone}(f)[-1],N) \cong \textrm{H}^{1}\,(\textrm{Hom}_R(\mathrm{cone}(f), N)) =0$.
Moreover, $\mathrm{Ext}_{\mathcal{X}\mathcal{M}}^{d+1}(M, N)$ vanishes by \cite[Lemma 4.5.(b)(2)]{16}.
Thus, we have the desired exact sequence.
\end{proof}

According to \cite[Fact 2.6]{18},
we know that $\mathcal {P}_C\perp\mathcal {P}_C$,
and $\mathcal {P}_C$ is exact and closed under kernels of epimorphisms.
Hence, it follows from Remark \ref{rem 3.12} that
$\mathcal {G}(\mathcal {P}_C)$ is also exact and closed under kernels of epimorphisms,
and $\mathcal {P}_C $ is both an injective cogenerator and a projective generator for $\mathcal {G}(\mathcal {P}_C)$.
Then we have the following result by Proposition \ref{prop 4.6}.

\begin{cor}\label{cor 4.7}{\cite[Theorem B]{18}}
Let $R$ be a commutative ring and $C$ a semidualizing $R$-module.
Assume that $M$ and $N$ are $R$-modules with $M\in\res\,\widehat{\mathcal {G}(\mathcal {P}_C)}$.
Set $d=\mathcal {G}(\mathcal {P}_C)$-$\pd(M)$.
Then there exists an exact sequence
$$0\longrightarrow \mathrm{Ext}_{\mathcal {G}(\mathcal {P}_C)\mathcal{M}}^{1}(M,
N)\longrightarrow\mathrm{Ext}_{\mathcal{P}_C\mathcal{M}}^{1}(M,
N)\longrightarrow\widehat{\mathrm{Ext}}_{\mathcal{P}_C\mathcal{M}}^{1}(M, N)\indent$$
$$\indent\longrightarrow\mathrm{Ext}^{2}_{\mathcal {G}(\mathcal {P}_C)\mathcal{M}}(M,
N)\longrightarrow\mathrm{Ext}_{\mathcal{P}_C\mathcal{M}}^{2}(M, N)
\longrightarrow\widehat{\mathrm{Ext}}_{\mathcal{P}_C\mathcal{M}}^{2}(M, N)\longrightarrow$$
$$\indent\cdots\longrightarrow\mathrm{Ext}^{d}_{\mathcal {G}(\mathcal {P}_C)\mathcal{M}}(M,
N)\longrightarrow\mathrm{Ext}_{\mathcal{P}_C\mathcal{M}}^{d}(M, N)\longrightarrow\widehat{\mathrm{Ext}}_{\mathcal{P}_C\mathcal{M}}^{d}(M, N)\longrightarrow0.
$$
\end{cor}

\begin{rem}\label{rem 2.1}
Under appropriate dual hypotheses, all the results in this paper have their dual versions.
\end{rem}

\bigskip \leftline {\bf Acknowledgments}\vspace{3mm}

This work is supported by the National Natural Science Foundation of China (11201063) and (11371089),
and the Specialized Research Fund for the Doctoral Program of Higher Education (20120092110020).

\end{document}